\documentclass[12pt]{article}
\usepackage{palatino} 
\newcommand{\reff}[1]{eq.~(\ref{#1})}

\begin{document}
\begin{titlepage}
\vskip.3in

\begin{center}
{\Large \bf Twisting invariance of link polynomials derived 
from ribbon quasi-Hopf algebras}
\vskip.3in
{\large J.R. Links, M.D. Gould and Y.-Z. Zhang} 
\vskip.2in
{\em Department of Mathematics, The University of Queensland, Brisbane,
     Qld 4072, Australia

Email: jrl@maths.uq.edu.au, yzz@maths.uq.edu.au}
\end{center}

\vskip 2cm
\begin{center}
{\bf Abstract}
\end{center}

The construction of link polynomials associated with finite dimensional
representations of ribbon quasi-Hopf algebra is discussed in terms of
the formulation of an appropriate Markov trace. We then show that this
Markov trace is invariant under twisting of the quasi-Hopf structure, 
which in turn implies twisting invariance of the associated link
polynomials.

\vskip 3cm

\end{titlepage}


\def\a{\alpha}
\def\b{\beta}
\def\d{\delta}
\def\e{\epsilon}
\def\ve{\varepsilon}
\def\g{\gamma}
\def\k{\kappa}
\def\l{\lambda}
\def\o{\omega}
\def\t{\theta}
\def\s{\sigma}
\def\D{\Delta}
\def\L{\Lambda}
\def\X{\bar{X}}
\def\Y{\bar{Y}}
\def\Z{\bar{Z}} 
\def\ch{\check}
\def\f{{\cal F}} 
\def\G{{\cal G}}
\def\hG{{\hat{\cal G}}}
\def\R{{\cal R}}
\def\hR{{\hat{\cal R}}}
\def\C{{\bf C}}
\def\P{{\bf P}}
\def\T{{\cal T}}
\def\H{{\cal H}}
\def\trho{{\tilde{\rho}}}
\def\tphi{{\tilde{\phi}}}
\def\tT{{\tilde{\cal T}}}
\def\uqsnh{{U_q[\widehat{sl(n|n)}]}}
\def\uqs1h{{U_q[\widehat{sl(1|1)}]}}
\def\ot{\otimes}


\def\beq{\begin{equation}}
\def\eeq{\end{equation}}
\def\bea{\begin{eqnarray}}
\def\eea{\end{eqnarray}}
\def\ba{\begin{array}}
\def\ea{\end{array}}
\def\no{\nonumber}
\def\lt{\left}
\def\rt{\right}
\newcommand{\bq}{\begin{quote}}
\newcommand{\eq}{\end{quote}}

\newtheorem{Theorem}{Theorem}
\newtheorem{Definition}{Definition}
\newtheorem{Proposition}{Proposition}
\newtheorem{Lemma}{Lemma}
\newtheorem{Corollary}{Corollary}
\newcommand{\proof}[1]{{\bf Proof. }
        #1\begin{flushright}$\Box$\end{flushright}}

\newcommand{\sect}[1]{\setcounter{equation}{0}\section{#1}}
\renewcommand{\theequation}{\thesection.\arabic{equation}}

\sect{Introduction\label{intro}}

The introduction of quantum algebras by Jimbo \cite{ji} and Drinfeld 
\cite{d1} lead to many remarkable developments in diverse areas of mathematical
physics. One such was in the field of knot theory whereby a connection
between the Yang-Baxter equation and the
braid group was quickly established. The quantum algebras, being
examples of quasi-triangular Hopf algebras, provide very systematic
means to find solutions of the Yang-Baxter equation which in turn gives
rise to representations of the braid group. Through a Markov trace
formulation defined on each braid group representation,
an invariant polynomial can then be computed for the knot
or link associated with the closure of the braid \cite{r,t,wda,zgb}. 
Extensions to accomodate the case of quantum superalgebras can be found
in \cite{lgz,z}.

Around the same time as the appearance of quantum algebras 
was Jones's discovery of a new polynomial invariant \cite{jo}, 
an evaluation of
which may be undertaken through the simplest quantum algebra
$U_q(sl(2))$ in its minimal
(two-dimensional) representation. 
After this breakthrough researchers  proceeded to obtain
generalizations with the notable examples being the HOMFLY \cite{h}
and Kauffman \cite{k} invariant polynomials.
What soon became apparent
was that the series of link polynomials associated with the fundamental
representations of the (non-exceptional) quantum algebras and superalgebras 
coincided with the two-variable
invariants developed in the wake of the discovery of Jones.  
More precisely,  the invariants associated with the fundamental
representations of the $U_q(gl(m|n))$ (which includes $U_q(gl(n))$) 
series belong to the class of HOMFLY
invariants while those of the $U_q(osp(m|2n))$ (including both
$U_q(o(m))$ and $U_q(sl(2n))$) series  
are of the Kauffman invariant type \cite{t,z}.    
It is important to emphasize, however, that by going to higher
representations new results are obtainable. In particular, the type I
quantum superalgebras consisting of $U_q(gl(m|n))$ and $U_q(osp(2|2n))$
admit one-parameter families of typical representations which give rise
to two-variable link invariants in a natural way \cite{dkl,glz,lg}. 
The work of Reshetikhin and Turaev \cite{rt} introduced further 
the notion of a
ribbon Hopf algebra as a particular type of quasi-triangular Hopf
algebra. All the quantum algebras fall into the class of ribbon Hopf
algebras. The algebraic properties of ribbon Hopf algebras is such that
an extension to produce invariants of oriented tangles is permissible.
A tangle diagram is analogous to a link diagram with the possibilty of
having free ends. An associated invariant takes the form of a tensor
operator acting on a product of vector spaces. 

As a generalization of Hopf algebras Drinfeld proposed the concept of 
quasi-Hopf algebras \cite{Dri90} whereby co-associativity of the co-algebra 
structure is not assumed. 
Any quasi-Hopf algebra generally belongs to an equivalence class where
each member is related to the others by twisting \cite{Dri90}. 
The potential for applications of these structures 
in the study of integrable systems is vast. 
They 
underly elliptic quantum algebras 
\cite{Enr97,Fel95,Fod94,Fro97,Jim97,Zha98} which play an important
role in obtaining solutions to the dynamical Yang-Baxter equations
\cite{Arn97,Bab96} and also twisting lies at the core of 
the construction of multiparametric
quantum spin chains \cite{flr}.

In the context of knot theory, Altschuler and Coste \cite{ac} 
have identified the corresponding ribbon quasi-Hopf algebras as the
necessary underlying algebraic structure to study tangle invariants 
(including closed link invariants). 
Here we wish to make two important observations to this field of study.
First, we will show that the class of ribbon quasi-Hopf algebras is
closed under twisting; i.e a twisted ribbon quasi-Hopf algebra is again
a ribbon quasi-Hopf algebra. Secondly, we assert that the link
polynomials computed from any finite dimensional representation of a
quasi-Hopf algebra are invariant with respect to twisting.
Importantly, this implies that link polynomials obtained from twisting
the usual quantum algebras give us nothing new. In this respect, one
cannot find twist generalizations of the HOMFLY and Kauffman invariants.
For a very special class of twists this result has already been noted by
Reshetikhin \cite{r1}, in which case the twisted quantum algebra is again
a ribbon Hopf algebra. Here we will prove the twisting invariance in
full generality.

The paper is structured as follows. We begin by presenting the
definition of a quasi-Hopf algebra. Next we show how representations of
the braid group are obtained from a representation of any quasi-Hopf
algebra. The third section deals with defining an appropriate Markov
trace on each braid group element which then affords a means to obtain a
link invariant. Finally, we demonstrate that the Markov trace is
invariant under any twisting. 

\sect{Quasi-Hopf Algebras}

Let us briefly recall the defining relations for 
quasi-Hopf algebras \cite{Dri90}. 

\begin{Definition}\label{quasi-bi}: 
A quasi-Hopf algebra is a  unital associative algebra $A$
over a field $K$  which is equipped with algebra homomorphisms $\e: 
A\rightarrow K$ (co-unit), $\D: A\rightarrow A\otimes A$ (co-product),
an invertible element $\Phi\in A\otimes A\otimes A$ 
(co-associator), an 
algebra anti-homomorphism $S: A\rightarrow A$ (anti-pode) and 
canonical elements $\a,~\b\in A$, satisfying 
\bea
&& (1\otimes\D)\D(a)=\Phi^{-1}(\D\otimes 1)\D(a)\Phi,~~
	\forall a\in A,\label{quasi-bi1}\\
&&(\D\otimes 1\otimes 1)\Phi \cdot (1\otimes 1\otimes\D)\Phi
	=(\Phi\otimes 1)\cdot(1\otimes\D\otimes 1)\Phi\cdot (1\otimes
	\Phi),\label{quasi-bi2} \label{pent}\\
&&(\e\otimes 1)\D=1=(1\otimes\e)\D,\label{quasi-bi3}\\
&&(1\otimes\e\otimes 1)\Phi=1,\label{quasi-bi4}\\
&& m\cdot (1\otimes\a)(S\otimes 1)\D(a)=\e(a)\a,~~~\forall
    a\in A,\label{quasi-hopf1}\\
&& m\cdot (1\otimes\b)(1\otimes S)\D(a)=\e(a)\b,~~~\forall a\in A,
     \label{quasi-hopf2}\\
&& m\cdot (m\otimes 1)\cdot (1\otimes\b\otimes\a)(1\otimes S\otimes
     1)\Phi^{-1}=1,\label{quasi-hopf3}\\
&& m\cdot(m\otimes 1)\cdot (S\otimes 1\otimes 1)(1\otimes\a\otimes
     \b)(1\otimes 1\otimes S)\Phi=1.\label{quasi-hopf4}
\eea
\end{Definition}
Here $m$ denotes the usual product map on $A$: $m\cdot (a\otimes b)=ab,~
\forall a,b\in A$. Note that since $A$ is associative we have
$m\cdot(m\otimes 1)=m\cdot (1\otimes m)$.
For all elements $a,b\in A$, the antipode satisfies
\beq
S(ab)=S(b)S(a).  
\eeq
The equations (\ref{quasi-bi2}), (\ref{quasi-bi3}) and (\ref{quasi-bi4}) imply
that $\Phi$ also obeys
\beq
(\e\otimes 1\otimes 1)\Phi=1=(1\otimes 1\otimes\e)\Phi.\label{e(phi)=1}
\eeq
Applying $\e$ to definition (\ref{quasi-hopf3}, \ref{quasi-hopf4}) 
we obtain, in view
of (\ref{quasi-bi4}), $\e(\a)\e(\b)=1$. 
By applying $\e$ to
(\ref{quasi-hopf1}), we have $\e(S(a))=\e(a),~\forall a\in A$.

A distinguishing feature of quasi-Hopf algebras is that they are in
general not co-associative; i.e
$$(1\otimes\D)\cdot\D\neq (\D\otimes 1)\cdot\D. $$
Thus for a given co-product the action extended to the $n$-fold tensor
product space is not uniquely determined. Throughout we will adopt the
convention to define a left co-product $\D_L$ which acts on the tensor
algebra $A^{\ot n}$ according to
$$\D_L(a\ot b\ot ....\ot c)=\D(a)\ot b\ot ....\ot c.  $$
We then recursively define the action
\beq \D^{(n)}=\D_L.\D^{(n-1)} \label{act} \eeq
with $\D^{(1)}=\D,\,\D^{(0)}=I$.

The category of quasi-Hopf algebras
is invariant under a kind of gauge transformation known as twisting. 
Let $(A,\D,\e,\Phi)$
be a quasi-Hopf algebra, with $\a,\b, S$ satisfying
(\ref{quasi-hopf1})-(\ref{quasi-hopf4}), and let $F\in A\otimes A$
be an invertible element satisfying the co-unit properties
\beq
(\e\otimes 1)F=1=(1\otimes \e)F.\label{e(f)=1}
\eeq
Throughout we set
\bea
&&\D_F(a)=F\D(a)F^{-1},~~~\forall a\in A,\label{twisted-d}\\
&&\Phi_F=F_{12}(\D\otimes
    1)F\cdot\Phi\cdot(1\otimes\D)F^{-1}F_{23}^{-1}\label{twisted-phi}
\eea 
where the subscripts above refer to the embedding of the elements in the
triple tensor product space. Then 
\begin{Theorem}\label{t-quasi-hopf}:
$(A,\D_F,\e,\Phi_F)$ defined by (\ref{twisted-d}, 
\ref{twisted-phi}) together with
$\a_F,\b_F, S_F$ given by
\beq
S_F=S,~~~\a_F=m\cdot(1\otimes\a)(S\otimes 1)F^{-1},~~~
	 \b_F=m\cdot(1\otimes\b)(1\otimes S)F,\label{twisted-s-ab}
\eeq
is also a quasi-Hopf algebra. Throughout, the element $F$ is referred to as
a  twistor.
\end{Theorem}

\begin{Definition}\label{quasi-quasi}: A quasi-Hopf
algebra $(A,\D,\e,\Phi)$ is called quasi-triangular if there
exists an invertible element $\R\in A\otimes A$ such that
\bea
&&\D^T(a)\R=\R\D(a),~~~~\forall a\in A,\label{dr=rd}\\
&&(\D\otimes 1)\R=\Phi^{-1}_{231}\R_{13}\Phi_{132}\R_{23}\Phi^{-1}_{123},
   \label{d1r}\\
&&(1\otimes \D)\R=\Phi_{312}\R_{13}\Phi^{-1}_{213}\R_{12}\Phi_{123}.
   \label{1dr}
\eea
We refer to $\R$ as the universal R-matrix.
\end{Definition}
Throughout, $\D^T=T\cdot\D$ with $T$ being the twist map
which is defined by 
\beq
T(a\otimes b)=b\otimes a;
\eeq
and $\Phi_{132}$ {\it etc} are derived from $\Phi\equiv\Phi_{123}$
with the help of $T$
\bea
&&\Phi_{132}=(1\otimes T)\Phi_{123},\no\\
&&\Phi_{312}=(T\otimes 1)\Phi_{132}=(T\otimes 1)
   (1\otimes T)\Phi_{123},\no\\
&&\Phi^{-1}_{231}=(1\otimes T)\Phi^{-1}_{213}=(1\otimes T)
   (T\otimes 1)\Phi^{-1}_{123},\no
\eea
and so on. 

It is easily shown that the properties (\ref{dr=rd})-(\ref{1dr})
imply the  Yang-Baxter type equation,
\beq
\R_{12}\Phi^{-1}_{231}\R_{13}\Phi_{132}\R_{23}\Phi^{-1}_{123}
  =\Phi^{-1}_{321}\R_{23}\Phi_{312}\R_{13}\Phi^{-1}_{213}\R_{12},
  \label{quasi-ybe}
\eeq
which is referred to as the  quasi-Yang-Baxter equation.  

\begin{Theorem}\label{t-quasi-quasi}: Denoting by the set
$(A,\D,\e,\Phi,\R)$  a
quasi-triangular quasi-Hopf algebra, then $(A, \D_F, \e, \Phi_F, \R_F)$
is also a quasi-triangular quasi-Hopf algebra, with the choice of
$R_F$ given by
\beq
\R_F=F^T \R F^{-1},\label{twisted-R}
\eeq
where $F^T=T\cdot F\equiv F_{21}$. Here $\D_F$ and $\Phi_F$ are given
by (\ref{twisted-d}) and (\ref{twisted-phi}), respectively.
\end{Theorem}

Let us specify some notations, where we adopt a summation convention
over all repeated indices. Throughtout the paper,
\bea
&&\Phi= X_i\otimes Y_i\otimes Z_i,~~~~
  \Phi^{-1}= \bar{X}_i\otimes\bar{Y}_i\otimes\bar{Z}_i,\no\\
&&F= f_i\otimes f^i,~~~~F^{-1}= \bar{f}_i\otimes \bar{f}^i,\no\\
&&\R=a_\nu\otimes b_\nu,~~~~\R^{-1}=c_\nu\otimes d_\nu,\no\\
&&(1\otimes\D)\D(a)=\sum a_{(1)}\otimes\D(a_{(2)})=\sum
   a^R_{(1)}\otimes a_{(2)}^R\otimes a^R_{(3)},\no\\
&&(\D\otimes 1)\D(a)=\sum \D(a_{(1)})\otimes a_{(2)})=\sum
   a^L_{(1)}\otimes a_{(2)}^L\otimes a^L_{(3)}, \no \\ 
&&(\Phi^{-1}\ot I).(\D\ot I\ot I)\Phi= A_i\ot B_i \ot C_i \ot D_i,
\no \\ 
&&(\D\ot I\ot I)\Phi^{-1}.(\Phi\ot I)= K_i\ot L_i\ot M_i\ot N_i
.\label{notation}
\eea

A important type of twistor is that due to Drinfeld \cite{Dri90}. For any
quasi-Hopf algebra $A$ observe that the actions
$$(S\ot S)\cdot\D^T,~~~~~~~~~~~~~~~\D^T\cdot S^{-1}$$
both determine algebra anti-homomorphisms. It follows that
$$\D'\equiv  (S\ot S)\cdot\D^T\cdot S^{-1}$$
gives rise to an algebra homomorphism and thus a co-product action on
$A$. Drinfeld showed that the actions $\D$ and $\D'$ are related by a
twistor; i.e.
$$\D(a)=\f^{-1}\left((S\ot S)\D^T(S^{-1}(a))\right) \f~~~~~~~~\forall
a\in A$$
where
$$ \f= (S\ot S)\D^T(X_i).\gamma.\D\left(Y_i\beta S(Z_i)\right) $$
and
\beq \gamma= S(B_i)\alpha C_i\ot S(A_i)\alpha D_i. \label{gamma} \eeq
It is also useful to define
\beq \delta= K_i\beta S(N_i)\ot L_i \beta S(M_i). \label{delta}  \eeq
Then the following relations can be shown to hold
$$\D(\alpha)=\f^{-1}\gamma,~~~~~~~~~\D(\beta)=\delta\f.$$

A quasi-Hopf algebra is said to be of trace type if
there exists an invertible    
element $u\in A$ such that
\beq
S^2(a)=u au^{-1},~~~~\forall a\in A.\label{s2a=u}
\eeq
In the case $A$ is quasi-triangular with R-matrix as in (\ref{notation})
we have \cite{ac} 
\begin{Theorem}\label{u-operator}: The operator defined by
\beq
u= S\lt(Y_i\b S(Z_i)\rt)S(\b_\nu)\a a_\nu X_i \label{u}
  \eeq
  satisfies (\ref{s2a=u}). Moreover the inverse is given by
  \beq
  u^{-1}= S^{-1}(X_i)S^{-1}(\a d_\nu)c_\nu 
    Y_i\b S(Z_i).\label{u-1}
    \eeq
    \end{Theorem}
An important relation satisfied by $u$ is 
\beq S(\alpha)u=S(b_\nu)\alpha a_\nu \label{imp} \eeq  
which we will need later.

The significance of trace type quasi-Hopf algebras is that they afford a
systematic means to construct Casimir invariants. We have the following
result from \cite{gzi}. 
\begin{Theorem}\label{trace-inv}: Let $\pi$ be the representation
afforded by the finite-dimensional  $A$-module $V$. Suppose
$\eta=\mu_i\otimes \nu_i\otimes \rho_i\in A\otimes {\rm End}V\otimes A$
obeys
\beq
(1\otimes\pi\otimes 1)(1\otimes\D)\D(a)\cdot\eta=\eta\cdot
 (1\otimes\pi\otimes 1)(1\otimes\D)\D(a),~~~~\forall a\in A,
 \label{eta}\eeq
 then
 \beq
 {\rm tr}\lt(\nu_i\pi\lt(\b S(\rho_i)S(\alpha)u\rt)\rt)\mu_i
 \eeq
 is a central element. Similarly if $\bar{\eta}=\bar{\mu}_i\otimes
 \bar{\nu}_i\otimes\bar{\rho}_i\in A\otimes {\rm End}V\otimes A$ satisfies
 \beq
 \bar{\eta}\cdot(1\otimes\pi\otimes 1)(\D\otimes 1)\D(a)=
   (1\otimes\pi\otimes 1)(\D\otimes 1)\D(a)\cdot\bar{\eta},~~~~\forall
      a\in A
      \eeq
      then
      \beq
      \sum {\rm
      tr}\lt(\bar{\nu}_i\pi\lt(u^{-1}S(\b)S(\bar{\mu}_i)\a\bar{\nu}_i\rt)
      \rt)\bar{\rho}_i
      \eeq
      is a central element.
      \end{Theorem}

As a consequence of the above we have 
\begin{Corollary}: Suppose $\o=\sum \o_i\otimes \Omega^i\in A\otimes{\rm
End}V$ satisfies
\beq
(1\otimes\pi)\D(a)\cdot\o=\o\cdot(1\otimes\pi)\D(a),~~~~\forall a\in A.
\eeq
Then (\ref{eta}) implies that
\bea
\tau(\o)&=&{\rm tr}\lt(\Omega^i\pi\lt(Y_k\b S(\bar{Z}_j
  Z_k) S(\a)u\bar{Y}_j  
  \rt)\rt)\bar{X}_j\o_i X_k\no\\
    \label{trace-o-inv}
    \eea
    is a central element.  
    
	  \end{Corollary}
For an $(n+1)$-fold tensor product space and $\o=\sum \o_i\otimes
\Omega^i\in A^{\ot n}\otimes{\rm
End}V$ we define 
\beq \tau_n(\o)={\rm tr}\lt(\Omega^i\pi\lt(Y_k\b S(\bar{Z}_j
  Z_k) S(\a)u\bar{Y}_j
    \rt)\rt)\D^{(n-1)}(\bar{X}_j)\o_i \D^{(n-1)}(X_k). \label{tn} \eeq

\sect{Representations of the braid group}

Given any representation $\pi$ of a
quasi-Hopf algebra $A$ we set 
\beq \check R=P.(\pi \ot \pi)\R \label{rh} \eeq  
and 
$$\Phi_i=(\D^{(i-2)}\ot I\ot I)\Phi.$$ 
In terms of $\check R$ the quasi-Yang-Baxter equation may be written
\beq \Phi\check R_{23}\Phi^{-1}\check R_{12}\Phi\check R_{23}\Phi^{-1}
=\check R_{12}\Phi\check R_{23}\Phi^{-1}\check R_{12} \label{qyb} \eeq  
where throughout we use the same symbols  $\Phi$ and $\Phi_i$ 
for both the algebraic objects and their matrix representatives.
\begin{Theorem} 
Define  $n$ operators on the $(n+1)$-fold tensor product space 
by 
\beq \sigma_i=\Phi_i\check R_{i(i+1)}\Phi_i^{-1},~~~~i=1,2,...,n\label{bg} 
\eeq    
These give
rise to a representation of the braid group $B_n$ 
by satisfying the defining relations
\bea \sigma_i\sigma_j&=&\sigma_j\sigma_i ~~~~~j\neq i\pm 1 \label{com} \\
\sigma_i\sigma_{i+1}\sigma_i&=&\sigma_{i+1}\sigma_i\sigma_{i+1}.\label{br}
\eea 
\end{Theorem}
The above result was given in \cite{ac}. Here we want to 
present a detailed proof.

First we establish that the braid generators (\ref{bg})
are invariant with respect to the co-product action
$\D^{(n)}$ of $A$; i.e
\beq [\sigma_i,\,\D^{(n)}(a)]=0~~~~~~~\forall a\in A.\label{inv} \eeq 
It is clear from the definition (\ref{rh}) that 
$$ [\check R,\,\D(a)]=0~~~~~~~\forall a\in A$$ 
which immediately implies that 
$$[\sigma_1,\,\D^{(n)}(a)]=0~~~~~~~\forall a\in A.$$ 
Next consider
\bea 
\s_2\D^{(j)}(a)&=&\Phi\check R_{23}\Phi^{-1}(\D\ot I^{\ot(j-1)})\D^{
(j-1)}(a) \no \\
&=&\Phi\check R_{23}(I\ot \D\ot I^{(j-2)})\D^{(j-1)}(a)\Phi^{-1} \no \\
&=&\Phi(I\ot \D\ot I^{(j-2)})\D^{(j-1)}(a)\check R_{23}\Phi^{-1} \no \\
&=&(\D\ot I^{\ot(j-1)})\D^{(j-1)}(a)\Phi\check R_{23}\Phi^{-1}\no \\ 
&=&\D^{(j)}(a) \s_2 \label{si2} \eea 
Observing that the action (\ref{act}) enjoys the property
$$\D^{(i)}\cdot \D^{(j)}=\D^{(i+j)}.$$
and applying $\D^{(k)}\ot I^{\ot j}$ to (\ref{si2}) now yields 
(\ref{inv}) by choosing $k=i-2,\,j=n-i-2$.

Since $\check R$ commutes with the co-product action we immediately
deduce for $i>1$   
\bea 
\s_1\s_i&=&\check R_{12} \Phi_i \check R_{i(i+1)} \Phi_i^{-1} \no \\
&=& \Phi_i\check R_{12}\check R_{i(i+1)}\Phi_i^{-1} \no \\ 
&=&\Phi_i\check R_{i(i+1)}\check R_{12}\Phi_i^{-1}\no \\ 
&=&\Phi_i\check R_{i(i+1)}\Phi_i^{-1}\check R_{12}\no \\ 
&=&\s_i\s_1. \no \eea 
Consider now for $l>3$   
\bea 
\s_2\s_l&=&\s_2 (\D^{(l-2)}\ot I\ot I)\Phi. \check R_{l(l+1)}(\D^{(l-2)}
\ot I \ot I)\Phi^{-1} 
\no \\
&=&(\D^{(l-2)}\ot I\ot I)\Phi.\s_2\check R_{l(l+1)}(\D^{(l-2)}\Phi^{-1}
\ot I\ot I)\no \\ 
&=&(\D^{(l-2)}\ot I\ot I)\Phi. \Phi_{123}\check R_{23}\Phi_{123}^{-1} 
\check R_{l(l+1)}(\D^{(l-2)}\ot I\ot I)\Phi^{-1} \no \\ 
&=&(\D^{(l-2)}\ot I\ot I)\Phi.\check R_{l(l+1)}\Phi_{123}
\check R_{23}\Phi_{123}^{-1}
(\D^{(l-2)}\ot I\ot I)\Phi^{-1} \no \\ 
&=&(\D^{(l-2)}\ot I\ot I)\Phi.\check R_{l(l+1)}\s_2(\D^{(l-2)} 
\ot I\ot I)\Phi^{-1} \no \\
&=& (\D^{(l-2)}\ot I\ot I)\Phi. \check R_{l(l+1)}(\D^{(l-2)} 
\ot I\ot I)\Phi^{-1}.\s_2 \no \\
&=& \s_2\s_l.   \label{s2com} \eea 
Applying  $\D^{(k)}\ot I\ot I$ to (\ref{s2com}) yields (\ref{com}) for
$i\geq 2$ by choosing $k=i-2,\,l=j-i+2.$  

In order to show that (\ref{br}) is satisfied we see from (\ref{qyb})
that 
$$\s_1\s_2\s_1=\s_2\s_1\s_2$$ 
is certainly true. Now through (\ref{qyb}), the invariance of 
$\check R$ and repeated use of the
pentagonal relation (\ref{pent}) we find 
\bea 
&&\s_2\s_3\s_2 \no   \\ 
&&=\Phi_2\ch R_{23}\Phi_2^{-1}\Phi_3\ch R_{34}\Phi_3^{-1}
\Phi_2\ch R_{23}\Phi^{-1}_2 \no \\
&&=\Phi_2\ch R_{23}\Phi_2^{-1}\Phi_3(I\ot I\ot\D)\Phi.\ch R_{34}
(I\ot I\ot\D)\Phi^{-1}.\Phi_3^{-1}
\Phi_2\ch R_{23}\Phi^{-1}_2 \no \\
&&=\Phi_2\ch R_{23}(I\ot \D\ot I)\Phi.(I\ot\Phi)
\ch R_{34}(I\ot\Phi^{-1})(I\ot \D\ot I)\Phi^{-1}
.\ch R_{23}\Phi^{-1}_2 \no \\
&&=\Phi_2(I\ot \D\ot I)\Phi.\left[\ch R_{23}(I\ot\Phi)
\ch R_{34}(I\ot\Phi^{-1}) \ch R_{23}\rt] 
(I\ot \D\ot I)\Phi^{-1}.\Phi^{-1}_2 \no \\  
&&=\Phi_2(I\ot \D\ot I)\Phi.\left[(I\ot\Phi)\ch R_{34}(I\ot\Phi^{-1}) 
\ch R_{23}(I\ot\Phi)\ch R_{34}(I\ot\Phi^{-1})\rt] 
(I\ot \D\ot I)\Phi^{-1}.\Phi^{-1}_2 \no \\  
&&=\Phi_3(I\ot I\ot \D)\Phi.\ch R_{34}(I\ot\Phi^{-1})
\ch R_{23}(I\ot\Phi)\ch R_{34}(I\ot I\ot \D)\Phi^{-1}.
\Phi_3^{-1} \no \\
&&=\Phi_3\ch R_{34}(I\ot I\ot \D)\Phi.(I\ot\Phi^{-1})
\ch R_{23}(I\ot\Phi)(I\ot I\ot \D)\Phi^{-1}.\ch R_{34}
\Phi_3^{-1}\no \\ 
&&=\Phi_3\ch R_{34}\Phi_3^{-1}\Phi_2
(I\ot \D\ot I)\Phi.\ch R_{23}(I\ot \D\ot I)\Phi^{-1}.\Phi^{-1}_2
\Phi_3\ch R_{23}\Phi_3^{-1}\no \\ 
&&=\Phi_3\ch R_{34}\Phi_3^{-1}\Phi_2\ch R_{23}\Phi_2^{-1}\Phi_3\ch
R_{34}\Phi_3^{-1} \no \\
&&=\s_3\s_2\s_2   
\label{int1} \eea  
Finally, acting $\D^{(i-2)}\ot I^{\ot 3}$ on (\ref{int1}) above yields 
$$\s_i\s_{i+1}\s_i=\s_{i+1}\s_i\s_{i+1}. $$

\sect{Link Polynomials from Ribbon Quasi-Hopf Algebras}

In \cite{ac} the following definition was proposed for the ribbon
quasi-Hopf algebras.
\begin{Definition} 
Let $A$ be a quasi-triangular quasi-Hopf algebra. We say that $A$ is a
ribbon quasi-Hopf algebra if there exists a central element $v\in A$
such that 
\begin{itemize}
\item[1.] $v^2=uS(u)$  
\item[2.] $S(v)=v$      
\item[3.] $\e (v)=1$ 
\item[4.] $\D(uv^{-1})=\f^{-1}(S\ot S)\f_{21}.uv^{-1}\ot uv^{-1}$ 
\end{itemize}
where $\f$ is the Drinfeld twist discussed earlier.
\end{Definition} 

Given a ribbon quasi-Hopf algebra, a prescription was also provided in
\cite{ac} to define a Markov trace on the braid group representation which 
in turn may be used to compute link polynomials in the usual 
way. From here on we will omit the symbol $\pi$ denoting the
representation for ease of notation.

\begin{Theorem} Let $\Psi$ be a word in the braid generators (\ref{bg})
for a fixed finite dimensional irreducible 
representation of a ribbon Hopf-algebra $A$. Then a Markov
trace $\theta_n $ on the $(n+1)$-fold tensor product space may be defined by 
$$\theta_n(\Psi)={\rm tr}\left(\Psi\D^{(n-1)}(\beta S(\alpha)uv^{-1}) 
\right)$$ 
which satisfies the Markov properties 
\begin{itemize}
\item[1.] $\theta_n(\Psi_1\Psi_2)=\theta_n(\Psi_2\Psi_1)
~~~~~\forall\,\Psi_1,\,\Psi_2\in B_n$ 
\item[2.] $\theta_n(\Psi\sigma^{\pm 1})=z^{\pm}\theta_{n-1}(\Psi) ~~~~~\forall
\Psi\in B_{n-1}\subset B_n$ 
\end{itemize}
where $z^{\pm}$ are the eigenvalues of the central operators $v^{\mp 1}$
in the representation $\pi$. 
\end{Theorem} 

The importance of the Markov trace is that from it one can 
define a link
polynomial $L(\hat{\Psi})$ through
\beq
L(\hat{\Psi})=(z^+z^-)^{n/2}\;\lt(\frac{z^-}{z^+}\rt)^{e(\t)/2}\t(\Psi),~~~~
  \t\in B_n   
  \eeq
  where $e(\Psi)$ is the sum of the exponents of the $\s_i$'s appearing in
  $\Psi$.
  The functional $L(\hat{\Psi})$ enjoys the following properties:
\begin{itemize}
\item[1.] $~~L(\widehat{\Psi\eta})=L(\widehat{\eta\Psi}),~\forall \Psi,
     \eta\in B_M;$ ~
\item[2.] $~~L(\widehat{\Psi\s_{n-1}^{\pm 1}})=L(\hat{\Psi}),~
	\forall \Psi\in B_{n-1}\subset B_n$
\end{itemize} 
and is an invariant of ambient isotopy.

The first Markov property follows easily from the invariance of the braid
generators $\sigma^{\pm 1}$ and the cyclic rule of traces. To establish the
second Markov property requires some work and was stated in \cite{ac}
without proof. Here we provide the details.

Before proceeding, we need to determine the co-product action of the
element $S(\alpha)uv^{-1}$. Using the Drinfeld twistor we find
\bea 
\D\left(S(\alpha)\rt)&=&\f^{-1}\lt((S\ot S)\D^T(\alpha)\rt)\f \no \\
&=&\f^{-1}\lt((S\ot S)(\f^{-1}_{21}\gamma_{21})\rt)\f \no \\
&=&\f^{-1}(S\ot S)\gamma_{21}.(S\ot S)\f_{21}^{-1}.\f \no \eea
Now through using (\ref{gamma}) and definition 3 we find that 
\beq \D\lt(S(\alpha)uv^{-1}\rt)= \f^{-1}\lt(S(D_i)S(\alpha)uv^{-1}A_i \ot
S(C_i)S(\a)uv^{-1}B_i\rt).  \label{dels}  \eeq 

We will also need the following result
\begin{Lemma} 
Let ${\cal C}\in {\rm End} (V\ot V)$ be any invariant operator; i.e 
$$[{\cal C}, \D(a)]=0 ~~~~~~\forall a\in A.$$ 
Then 
$$\lt(A_i\ot B_i\rt){\cal C}\lt(K_j\b S(D_iN_j)\ot L_j\b S(C_iM_j)\rt)
=\lt(\X_j\ot\Y_j\rt){\cal C}\lt(X_i\b\ot Y_i\b S(\Z_j Z_i)\rt).$$
\end{Lemma}
The above result follows directly from the definitions (\ref{notation}). 
We may now see that 
\bea 
\theta_n(\Psi)&=&{\rm tr}\left(\Psi\D^{(n)}(\beta S(\alpha)uv^{-1})
\right)  \no \\
&=&{\rm tr}\left(\Psi.(\D^{(n-1)}\ot
I)\d.\D^{(n-1)}(S(D_i)S(\alpha)uv^{-1}A_i)\ot S(C_i)S(\alpha)
uv^{-1}B_i \right) \no \\
&=&{\rm tr}\left(\Psi.\D^{(n-1)}\left(K_j\beta S(N_j)S(D_i)S(\alpha)
uv^{-1}A_i\right) \ot
L_j\beta S(M_j)S(C_i)S(\alpha) uv^{-1}B_i  
\right)  \no \\ 
&=&{\rm tr}\left(\Psi. \D^{(n-1)}(X_i\beta S(\alpha)uv^{-1}\X_j)\ot
Y_i\beta S(Z_i)S(\Z_j)S(\alpha)uv^{-1}\Y_j  \rt) \no \\
&=&{\rm tr}\left(\tau_n(\Psi)\D^{(n-1)}(\beta S(\alpha)uv^{-1})\rt) \no \\
&=&\theta_{n-1}\left(\tau_n(\Psi)\right) \eea  
where the element $u$ in the definition (\ref{tn}) has now been replaced
by $uv^{-1}$. 
It is apparent also from (\ref{tn}) that for $\Psi\in B_{n-1}$
then 
$$\tau_n(\Psi\sigma_n^{\pm 1})=\Psi\tau_n(\sigma_n^{\pm 1})$$  
To evaluate $\tau_n(\sigma_n^{\pm 1})$ we can appeal to the pentagonal
relation (\ref{pent}) to find that 
\bea &&\Phi^{-1}_{n+1}\sigma_n^{\pm 1}\Phi_{n+1}\no \\ 
&&= 
\Phi^{-1}_{n+1}\Phi_n\ch R_{n(n+1)}^{\pm 1}\Phi^{-1}_n\Phi_{n+1} \no \\
&&=I^{(n-2)}\ot \left(( I \ot I \ot \D)\Phi.(I\ot \Phi^{-1})(I\ot \D\ot
I)\Phi^{-1}.(I\ot \ch R^{\pm 1} \ot I)\rt. \no \\ 
&&~~~~~~~~~~~.\lt.(I\ot \D\ot I)\Phi.(I\ot \Phi)
(I\ot I\ot \D)\Phi^{-1}\rt) \no \\  
&&=I^{(n-2)}\ot \left((I\ot I \ot \D)\Phi.(I\ot \Phi^{-1})
(I\ot \ch R^{\pm 1} \ot I)(I\ot \Phi)
(I\ot I\ot \D)\Phi^{-1}\rt) \no \eea 
which, upon using (\ref{quasi-hopf1}, \ref{quasi-hopf2}, 
\ref{e(phi)=1}),  leads us to conclude that 
$$\tau_n(\sigma_n^{\pm 1})=I^{\ot (n-1)}\ot \tau(\ch R^{\pm 1}
). $$ 
An algebraic exercise shows that 
\bea \tau(\ch R)&=&\X_jb_{\nu}Y_l\beta S(Z_l)S(\Z_j)S(\alpha)uv^{-1}\Y_j
a_{\nu}X_l, \no  \\
\tau(\ch R^{-1})&=&X_jc_{\nu}Y_k\beta
S(Z_k)S(\Z_j)S(\alpha)uv^{-1}\Y_jd_{\nu}X_k 
\no \eea  
and hence we can conclude that $z^{\pm}$ are given by the eigenvalues of
the central operators $v^{\mp 1}$ if we can show that the 
following relations hold.  
\begin{Lemma} 
\bea \X_jb_{\nu}Y_l\beta S(Z_l)S(\Z_j)S(\alpha)u\Y_j  
a_{\nu}X_l&=&I, \no   \\
\X_jc_{\nu}Y_k\beta S(Z_k)S(\Z_j)S(\alpha)u\Y_jd_{\nu}X_k 
&=&v^2.\no  \eea  
\end{Lemma} 
Through use of (\ref{quasi-hopf3}, \ref{imp}) we obtain   
\bea I&=&\X_i\beta S(\Y_i)\alpha \Z_i \no \\
&=&\X_j\beta S(\Y_i)S(d_{\nu})S(b_{\mu})\alpha a_{\mu}c_{\nu}\Z_i \no \\
&=&\X_j\beta S(\Y_i)S(d_{\nu})S(\alpha)uc_{\nu}\Z_i \no \\
&=&\X_i\beta S(\Y_i)S(d_{\nu})S(\alpha)S^2(c_\nu)S^2(\Z_i)u. 
\eea 
From \reff{1dr} we see that 
$$\R^{-1}_{13}\Phi^{-1}_{312}\left(I\ot \D\right)\R=\Phi^{-1}_{213}\R_{12}
\Phi_{123} $$ 
which expressed in terms of the tensor components reads 
$$c_{\nu}\Z_ja_l\ot \X_jb_l^{(1)}\otimes d_{\nu}\Y_j b_l^{(2)} 
=\Y_j a_\nu X_l\ot \X_j b_\nu Y_l\ot \Z_jZ_l. $$ 
We can now write 
\bea 
&&\X_jb_{\nu}Y_l\beta S(Z_l)S(\Z_j)S(\alpha)u\Y_j  
a_{\nu}X_l \no \\
&&~~~~=\X_jb_l^{(1)}\beta S(b_l^{(2)})S(\Y_j)S(d_\nu)S(\alpha)u
c_\nu\Z_ja_l \no \\
&&~~~~=\e(b_l)\X_j\beta S(\Y_j)S(d_\nu)S(\alpha)u
c_\nu\Z_ja_l \no \\  
&&~~~~=\X_i\beta S(\Y_i)S(d_{\nu})S(\alpha)S^2(c_\nu)S^2(\Z_i)u \no \\
&&~~~~=I.   \eea 

Next we see that 
\bea 
u&=& S\left(\X_i\beta S(\Y_i)\alpha\Z_i\right)u \no \\
&=&S(\Z_i)S(\alpha)S^2(\Y_i)S(\beta)S(\X_i)u \no \\
&=&S(\Z_i)S(\alpha)u\Y_iS^{-1}(\beta)S^{-1}(\X_i) \no \\
&=&S(\Z_i)S(b_\nu)\alpha a_\nu\Y_iS^{-1}(\beta)S^{-1}(\X_i) 
\eea 
where in the last step we have used \reff{imp}. Consequently
$$S(u)=\X_i\beta S(\Y_i)S(a_\nu)S(\alpha)S^2(b_\nu)S^2(\Z_i). $$ \
From \reff{d1r} we have 
$$\R_{23}\Phi^{-1}_{123}\left(\D\ot I\right)\R^{-1}=\Phi_{132}^{-1}
\R^{-1}_{13}\Phi_{231} $$ 
which we may express as 
$$\X_jc^{(1)}_\nu\ot a_\mu\Y_jc_\nu^{(2)}\ot b_\mu \Z_jd_\nu 
=\X_jc_\nu Y_k\ot \Z_jZ_k\ot \Y_jd_\nu X_k. $$ 
This relation leads us to deduce that 
\bea 
&&\X_jc_{\nu}Y_k\beta S(Z_k)S(\Z_j)S(\alpha)u\Y_jd_{\nu}X_k  
\no \\
&&~~~~=\X_jc_\nu^{(1)}\beta S(c_\nu^{(2)})S(\Y_j)S(a_\mu) S(\alpha)u
b_\mu\Z_jd_\nu \no \\
&&~~~~=\X_j\beta S(\Y_j)S(\a_\mu)S(\alpha)S^2(b_\mu)S^2(\Z_j)u \no \\
&&~~~~=S(u)u \no \\
&&~~~~=v^2 \eea 
which proves lemma 2 and completes the proof of Theorem 6. 

\sect{Twisting Invariance of the Markov trace}

Now we are in a position to show twisting invariance of the link
polynomials. Let us begin with the following result.
\begin{Proposition} 
Every twisted ribbon quasi-Hopf algebra is
again a ribbon quasi-Hopf algebra. 
\end{Proposition} 
Recall from definition 3 that the first three conditios of a ribbon
quasi-Hopf algebra are properties of the algebra structure rather than
the co-algebra. Thus, to this end we need only show that if
$$
 \Delta(uv^{-1})=\f^{-1}(S\ot S)\f_{21}.\,uv^{-1}\ot uv^{-1} $$ 
 then 
$$\Delta_F(uv^{-1})=\f_F^{-1}(S\ot S)(\f_F)_{21}.\,uv^{-1}\ot uv^{-1}$$ 
where  $\f_F$ denotes the Drinfeld twistor for the twisted quasi-Hopf
algebra. 
Recalling that the Drinfeld twist $\f$ is determined by 
$$\f\Delta(a)\f^{-1}=(S\ot S)\left(\Delta^T\left(S^{-1}(a)\right)\right)
~~~~~\forall a\in A $$ 
shows that 
\bea F\Delta(a)F^{-1} &=& \Delta_F(a)\no \\
&=&\f^{-1}_F\left((S\ot
S)\Delta^T_F\left(S^{-1}(a)\right)\right)\f_F \no \\
&=&\f_F^{-1}\left((S\ot S)\left(F_{21}\Delta^T\left(S^{-1}(a)\right)
F_{21}^{-1}\right)\right)\f_F \no \\
&=&\f_F^{-1}\left((S\ot S)F_{21}^{-1}.(S\ot
S)\Delta^T\left(S^{-1}(a)\right).(S\ot S)F_{21}\right) \f_F \no \\
&=&\f^{-1}_F(S\ot S)F_{21}^{-1}.\f\Delta(a)\f^{-1}.(S\ot S)F_{21}\f_F \no 
\eea 
which leads us to  
$$\f_F=(S\ot S) F_{21}^{-1}.\f F^{-1}. $$ 

Now we observe that 
\bea &&\f_F^{-1}(S \ot S)(\f_F)_{21}.uv^{-1}\ot uv^{-1} \no \\ 
&&~~~~~=F\f^{-1}(S\ot S)F_{21}.(S\ot S)\left((S\ot S)F^{-1}.(\f F^{-1})_{21}
\right)uv^{-1}\ot uv^{-1} \no \\  
&&~~~~~=F\f^{-1}(S\ot S)\f_{21}.(S^2\ot S^2)F^{-1}.uv^{-1}\ot uv^{-1} \no \\
&&~~~~~=F\f^{-1}(S\ot S)\f_{21}.uv^{-1}\ot uv^{-1}. F^{-1}\no \\
&&~~~~~=F\Delta(uv^{-1})F^{-1} \no \\
&&~~~~~=\Delta_F(uv^{-1}) \no \eea 
thus establishing that the twisted ribbon quasi-Hopf algebra is also of ribbon
type.

By induction the co-product action on the $(n+1)$-fold space assumes the form
$$\Delta^{(n)}_F(a)=\chi_{n}\Delta^{(n)}\chi_{n}^{-1}$$ 
where  
$$\chi_{n}=F_{12}.(\Delta\ot I)F_{12}.(\Delta^2\ot
I)F_{12}....(\Delta^{(n-1)}\ot I )F. $$ 
Consider next  
\bea 
\chi_{n}\sigma_i\chi^{-1}_{n}&=&F_{12}(\Delta\ot
I)F_{12}.(\Delta^2\ot I)F_{12}....(\Delta^{(n-1)}\ot I )F
\sigma_i \no \\ 
&& ~~~~~~~~~~  
(\Delta^{(n-1)}\ot I )F^{-1}....(\Delta\ot I)F_{12}^{-1}.F_{12}^{_1} \no
\\  
&=&F_{12}(\Delta\ot
I)F_{12}.(\Delta^2\ot I)F_{12}....(\Delta^{(i-1)}\ot I )F 
\sigma_i
\no \\ &&  ~~~~~~~~~ 
(\Delta^{(i-1)}\ot I )F^{-1}....(\Delta\ot I)F_{12}^{-1}.F_{12}^{_1} \no
\\
 &=&\chi_i\sigma_i\chi_i^{-1}. \no \eea  
   
We now determine the representations of the braid generators under
twisting; i.e. 
\bea 
\sigma_i^F&=& (\Delta_F^{(i-2)}\ot I\ot I)\Phi_F.\check R^F_{i(i+1)}.
(\Delta_F^{(i-2)}\ot I\ot I) \Phi_F^{-1} \no \\
&=&\chi_{i-2}(\Delta^{(i-2)}\ot I \ot I)\left(F_{12}(\Delta\ot I)F.\Phi.
(I\ot \Delta)F^{-1}.F_{23}^{-1}\right)\chi_{i-2}^{-1}F_{i(i+1)}\check
R_{i(i+1)}F^{-1}_{i.i+1} \no \\  
&& ~~~~~~~~~\chi_{i-2}\Delta^{(i-2)}\left(F_{23}(I\ot
\Delta)F.\Phi^{-1}(\Delta\ot I)F^{-1}.F_{12}^{-1}\right) \chi^{-1}_{i-2}
\no \\ 
&=&\chi_i\Delta\Phi.(\Delta^{(i-2)}\otimes\Delta)F^{-1}.F_{i(i+1)}^{-1}
\chi_{i-2}^{-1}F_{i(i+1)}\check R_{i(i+1)}F^{-1}_{i(i+1)} \no \\
&& ~~~~~~~~~\chi_{i-2} F_{i(i+1)}(\Delta^{(i-2)}\ot
\Delta)F.\Delta^{(i-2)}\Phi^{-1}. \chi_i^{-1} \no \\
&=&\chi_i\Delta^{(i-2)}\Phi. \check R_{i(i+1)} \Delta^{(i-2)} \Phi^{-1}.
\chi_i^{-1} \no \\ 
&=&\chi_i\sigma_i\chi_i^{-1} \no \\
&=&\chi_{n}\sigma_i\chi_{n}^{-1} \no \eea 
which shows that the representation of the braid generators under
twisting are related to those of the untwisted case by a basis
transformation. Thus for any word in the generators of the braid group
we can write 
$$\Psi^F=\chi_{n}\Psi\chi_{n}^{-1}$$ 
in an obvious notation. 

Using the relations (\ref{twisted-s-ab}) we may write 
$$\alpha_F=S(\bar{f}_i)\alpha \bar{f}^i,~~~~~~~~\beta_F=f_i\beta S(f^i)$$
and proceed to calculate 
\bea 
\theta^F_n(\Psi)&=&{\rm tr}\left(\Psi^F \Delta_F^{(n)}\left(\beta_F S(\alpha_F)
uv^{-1}\right)\right) \no \\ 
&=&{\rm tr}\left(\chi_{n}\Psi \Delta^{(n)}
\left(\beta_F S(\alpha_F)uv^{-1}\right)\chi^{-1}_{n}\right) \no \\
&=&{\rm tr}\left(\Psi\Delta^{(n)}\left(f_i\beta S(f^i)S(\bar{f}^j)S(
\alpha)S^2(\bar{f}_j)uv^{-1}\right)\right) \no \\
&=&{\rm tr}\left(\Psi\Delta^{(n)}\left(f_i\beta S(\bar{f}^jf^i)
S(\alpha)uv^{-1}\bar{f}_j\right)\right) \no \\
&=&{\rm tr}\left(\Delta^{(n)}(\bar{f}_j)\Psi\Delta^{(n)}\left(f_i\beta
S(\bar{f}^jf^i)S(\alpha)uv^{-1}\right)\right) \no \\
&=&{\rm tr}\left(\Psi\Delta^{(n)}\left(\bar{f}_jf_i\beta 
S(\bar{f}^jf^i)S(\alpha)uv^{-1}\right)\right) \no \\
&=&{\rm tr}\left(\Psi \Delta^{(n)}\left(\beta S(\alpha) uv^{-1} 
\right) \right) \no \\ 
&=&\theta_n(\Psi) \no \eea 
which proves twisting invariance of the Markov trace 
and consequently the associated link polynomials.

\vskip.3in
\noindent {\bf Acknowledgements.}

This work has been financed by the Australian Research Council.

\newpage

\end{document}